\relax
\documentclass[runningheads]{llncs}
\usepackage{times} 
\usepackage{helvet} 
\usepackage{courier} 
\usepackage[hyphens]{url} 
\usepackage{graphicx} 
\usepackage{graphicx}  
\usepackage{caption} 
\usepackage{amsfonts}
\usepackage{amsmath}
\usepackage{todonotes}
\usepackage{algorithm}
\usepackage{algorithmic}
\usepackage{subcaption}
\usepackage{float}

\begin{document}

\title{Continuous surrogate-based optimization algorithms are well-suited for expensive discrete problems}

\author{Rickard Karlsson \and Laurens Bliek \and Sicco Verwer \and Mathijs de Weerdt
}
\institute{Delft University of Technology,
              Faculty of Electrical Engineering, Mathematics and Computer Science, 
              Van Mourik Broekmanweg 6, 
              2628 XE Delft,
              The Netherlands}
\maketitle
\begin{abstract}
One method to solve expensive black-box optimization problems is to use a surrogate model that approximates the objective based on previous observed evaluations. The surrogate, which is cheaper to evaluate, is optimized instead to find an approximate solution to the original problem. In the case of discrete problems, recent research has revolved around surrogate models that are specifically constructed to deal with discrete structures. A main motivation is that literature considers continuous methods, such as Bayesian optimization with Gaussian processes as the surrogate, to be sub-optimal (especially in higher dimensions) because they ignore the discrete structure by, e.g., rounding off real-valued solutions to integers. However, we claim that this is not true. In fact, we present empirical evidence showing that the use of continuous surrogate models displays competitive performance on a set of high-dimensional discrete benchmark problems, including a real-life application, against state-of-the-art discrete surrogate-based methods. Our experiments on different discrete structures and time constraints also give more insight into which algorithms work well on which type of problem.

\keywords{Surrogate models, Bayesian optimization, expensive combinatorial optimization, black-box optimization}
\end{abstract}

\section{Introduction}

A principal challenge in optimization is to deal with black-box objective functions. The objective function is assumed to be unknown in this case, in contrast to traditional optimization that often utilizes an explicit formulation to compute the gradient or lower bounds. Instead, we assume to have an objective $y = f(\textbf{x}) + \epsilon$ with some unknown function $f(\textbf{x})$ together with additive noise $\epsilon$. Furthermore, $f(\textbf{x})$ can be expensive to evaluate in terms of time or another resource which restricts the number of evaluations allowed.

One type of method to solve these black-box optimization problems is the use of surrogate models. Surrogate-based algorithms approximate the objective function in search of the optimal solution, with the benefit that the surrogate model is cheaper to evaluate.
Bayesian optimization~\cite{snoek2012practical} is an example of such a surrogate-based algorithm.

An active field of research is how to deal with discrete black-box optimization problems with an expensive objective function. There are many real-world examples of this, such as deciding on the architecture of a deep neural network~\cite{elsken2018neural} or designing molecules with desirable properties~\cite{korovina2020chembo}. Furthermore, optimization over structured domains was highlighted as an important problem to address from the NIPS 2017 workshop on Bayesian optimization~\cite{hernandez-lobato2017NIPSworkshop}.

Discrete optimization problems can be solved with a continuous surrogate model, e.g., Bayesian optimization with Gaussian processes~\cite{snoek2012practical}, by ignoring the discrete structure and rounding off the real-valued input to discrete values. However, literature in this field  generally considers this to be a sub-optimal approach \cite{baptista2018bayesian,garrido2020dealing}. Therefore, research has revolved around inherently discrete models such as density estimators or decision trees, e.g., HyperOpt~\cite{bergstra2013hyperopt} or SMAC~\cite{hutter2011sequential}. Another approach is to use continuous models that guarantee discrete optimal solutions, such as the piece-wise linear model IDONE~\cite{bliek2020black}.

In contrast to common belief, we present an empirical study that displays that continuous surrogate models, in this case Gaussian processes and linear combinations of rectified linear units, show competitive performance on expensive discrete optimization benchmarks by outperforming discrete state-of-the-art algorithms. 
Firstly, we will introduce the problem, the related work, and the considered benchmark problems. Then, in the remainder of the paper we
1) perform a benchmark comparison  between continuous and discrete surrogate-based algorithms on optimization problems with different discrete structures (including one real-life application),
2) investigate why continuous surrogate models perform well by transforming the different discrete problem structures and visualizing the continuous surrogate models,
and 
3) perform a more realistic analysis that takes the time budget and evaluation time into account when comparing the algorithms. 
We conclude that continuous surrogates applied to discrete problems should get more attention, and leave some questions for interesting directions of future research in the domain of discrete expensive black-box optimization.

In this paper, we first present the problem description and surrogate-based optimization in Section~2. Then, Section~3 gives an overview of the related work and surrogate-based algorithms for discrete optimization problems. In Section~4 we explain the experimental setup and the combinatorial problems that we use in our benchmark comparison, and Section~5 presents the results from our experiments. Lastly, in Section~6 we give the conclusion from this work and propose future work. 

\section{Problem Description}

Consider the following class of \textit{d}-dimensional discrete optimization problems:

\begin{equation}
    \begin{aligned}
    & \underset{\mathbf{x}}{\text{minimize}} && f(\mathbf{x}) \\
    & \text{subject to } && \mathbf{x}\in\mathbb{Z}^d \\
    & && l_i \leq \mathbf{x}_i \leq u_i, \; i =1,\dots,d
    \end{aligned}
    \label{eq:mainproblem}
\end{equation}

where $l_i$ and $u_i$ are the lower and upper bound for each integer-valued decision variable $x_i$. 
For black-box optimization problems, we assume to have no closed form expression for $f:\mathbb{Z}^d\rightarrow\mathbb{R}$. 
The only information which can be gathered about $f$ comes from observing the output when evaluating $f(\mathbf{x})$ given some input $\mathbf{x}$. 
However, in many real-world applications we will also have to deal with some noise $\epsilon\in\mathbb{R}$ such that we are given the output $y = f(\mathbf{x}) + \epsilon$. 
Obtaining an evaluation is also assumed to be expensive: it could require large computational power, human interaction with the system or time consuming simulations.
Therefore it is of interest to obtain a solution within a limited amount of evaluations $B$, also known as the budget.

One way of solving this class of problems is to make use of a so called surrogate model.
A surrogate model is an auxiliary function $M$ that approximates the objective function based on the points evaluated so far.
This model is cheaper to evaluate in comparison to the original black-box objective function.
Given a number $m$ of already evaluated points, the surrogate model is constructed using the evaluation history $H~=~\{(\mathbf{x}^{(1)},y^{(1)}), (\mathbf{x}^{(2)},y^{(2)}), \dots, (\mathbf{x}^{(m)},y^{(m)})\}$.
The surrogate can be utilized to predict promising points to evaluate next on. The next feasible solution $\mathbf{x}^{(m+1)}$ to evaluate on can be chosen based on this prediction.
These steps, which are also described in Algorithm~\ref{alg:surrogate_modelling}, are repeated until the budget $B$ is reached.

Typically, an acquisition function $A(M,\mathbf{x})$ is used to propose the next point $\mathbf{x}^{(m+1)}$ to evaluate with the objective function. 
It predicts how promising a new point $\mathbf x$ is, based on a trade-off of exploitation (searching at or near already evaluated points that had a low objective) and exploration (searching in regions where the surrogate has high uncertainty).
In general, the next point is chosen by finding the global optimum $\mathbf{x}^{(m+1)}=\operatornamewithlimits{argmax}\limits_{\mathbf{x}} A(M,\mathbf{x})$.

\begin{algorithm}
\caption{Surrogate-based optimization}
\label{alg:surrogate_modelling}

    \begin{algorithmic}[1]
    \REQUIRE budget $B$, surrogate model $M$, acquisition function $A$
    \STATE Initialize $\mathbf{x}^{(1)}$ randomly and an empty set $H$
    \FOR{$m=1:B$}
    
        \STATE $y^{(m)} \gets f(\mathbf{x}^{(m)}) + \epsilon$
        \STATE $H\gets H \cup \{(\mathbf{x}^{(m)}, y^{(m)}))\}$
        \STATE $M \gets$ fit surrogate model using H
        \STATE $\mathbf{x}^{(m+1)} \gets \operatornamewithlimits{argmax}\limits_{\mathbf{x}} A(M,\mathbf{x})$
    \ENDFOR
    \RETURN optimal $(\mathbf{x}^*,y^*)\in H$
    
    \end{algorithmic}

\end{algorithm}
\section{Related Work}
\label{seq:related_work}

Although discrete problem structures are difficult to handle in black-box optimization, multiple approaches have been proposed. 
A survey by M. Zaefferer~\cite{zaefferer2018Auth} presents different strategies for dealing with discrete structures in surrogate-based algorithms. The first strategy is the naive way by simply ignoring the discrete structure.
Another strategy is to use inherently discrete models such as tree-based models~\cite{bergstra2013hyperopt,hutter2011sequential}. These models can however fail if the problem structure is too complex or if there are both discrete and continuous variables involved~\cite{zaefferer2018Auth}. Lastly, discrete structures can be dealt with by using a certain mapping. Although this strategy does not apply directly to a surrogate model, a suitable mapping can make the problem easier. For example, encoding integer solutions with a binary representation can be easier for some regression models to handle. 

There are also other strategies such as using problem-specific feature extraction or customizing the model. However, these violate the black-box assumption which is why we will not discuss them.

We now discuss several surrogate-based optimization algorithms that can solve the expensive discrete optimization problem in eq.~\eqref{eq:mainproblem} and that also have their code available online.

Bayesian optimization has a long history of success in expensive optimization problems~\cite{jones1998efficient}, and has been applied in many domains such as chemical design and hyperparameter optimization for deep learning~\cite{griffiths2017constrained,klein2017fast}.
It typically uses a Gaussian process as a surrogate to approximate the expensive objective.
Several acquisition functions exist to guide the search, such as Expected Improvement, Upper Confidence Bound, or Thompson sampling~\cite{shahriari2016reviewBO},
information-theoretic approaches such as Predictive Entropy Search~\cite{hernandez2014predictive}, or simply the surrogate itself~\cite{de2019greed,rehbach2020expected}.
Though Gaussian processes are typically used on continuous problems, they can be adapted for problems with discrete variables as well.
The authors of \cite{garrido2020dealing} suggest three possible approaches, namely rounding to the nearest integer 1) when choosing where to evaluate the objective function, 2) when evaluating the objective function, or 3) inside the covariance function of the Gaussian process.
The latter provides the best results but gives an acquisition function that is hard to optimize.
The first option leads to the algorithm getting stuck by repeatedly evaluating the same points, although this can be circumvented by carefully balancing exploration and exploitation~\cite{luong2019bayesian}.
In this work, we will consider only the simpler second option, for which we do not need to modify any existing implementations.\footnote{We consider the implementation from \url{https://github.com/fmfn/BayesianOptimization} in this work, which uses the Upper Confidence Bound acquisition function.}

BOCS\footnote{\url{https://github.com/baptistar/BOCS}}~\cite{baptista2018bayesian} transforms the combinatorial problem into one that can be solved with semi-definite programming.
It uses Thompson sampling as the acquisition function.
However, it suffers from a large time complexity, which was only recently alleviated by using a submodular relaxation called the PSR method\footnote{\url{https://github.com/aryandeshwal/Submodular_Relaxation_BOCS}}~\cite{deshwal2020scalable}.

COMBO\footnote{\url{https://github.com/tsudalab/combo}}~\cite{ueno2016combo} uses an efficient approximation of a Gaussian process with random features, together with Thompson sampling as the acquisition function. Though this gives increased efficiency, COMBO deals with discrete search spaces by iterating over all possible candidate solutions, which is only possible for small-dimensional problems.
Later, a different group proposed another algorithm with the same name\footnote{\url{https://github.com/QUVA-Lab/COMBO}}, based on the graph Fourier transform~\cite{oh2019combinatorial}.
However, this method uses approximately the same computational resources as BOCS.

HyperOpt\footnote{\url{https://github.com/hyperopt/hyperopt}}~\cite{bergstra2013hyperopt}
makes use of a tree of Parzen estimators as the surrogate model.
It can naturally deal with categorical or integer variables, and even with conditional variables that only exist if other variables take on certain values.
The algorithm is known to perform especially well on hyperparameter tuning problems with hundreds of dimensions~\cite{bergstra2013making}.
This is in sharp contrast with Bayesian optimization algorithms using Gaussian processes, which are commonly used on problems with less than $10$ dimensions.
A possible drawback for HyperOpt is that each dimension is modeled separately, i.e., no interaction between different variables is modeled.
HyperOpt uses the Expected Improvement acquisition function.

SMAC\footnote{\url{https://github.com/automl/SMAC3}}~\cite{hutter2011sequential} is another surrogate-based algorithm that can naturally deal with integer variables.
The main reason for this is that the surrogate model used in this algorithm is a random forest, which is an inherently discrete model.
A point of critique for SMAC is that the random forests have worse predictive capabilities than Gaussian processes.
Nevertheless, like HyperOpt, SMAC has been applied to problems with hundreds of dimensions~\cite{lindauer2017warmstarting}.
SMAC uses the Expected Improvement acquisition function.

IDONE\footnote{\url{https://bitbucket.org/lbliek2/idone}}~\cite{bliek2020black} uses a linear combination of rectified linear units as its surrogate model. 
This is a continuous function, yet it has the special property that any local minimum of the model is located in a point where all variables take on integer values.
This makes the method suitable for expensive discrete optimization problems, with the advantage that the acquisition function can be optimized efficiently with continuous solvers.
IDONE uses the surrogate model itself as the acquisition function, but adds small perturbations to the optimum of the acquisition function to improve its exploration capabilities.
Though the method is not as mature as SMAC or HyperOpt, it also has been applied to problems with more than $100$ variables~\cite{bliek2020black}.

\section{Benchmark problems}
\label{seq:experiments}
We present the four different benchmark problems that are used to compare the surrogate-based algorithms. The purpose of the benchmarks is to compare the discrete surrogate-based algorithms presented in the previous section and investigate which algorithms are most suited for which type of problem.

The benchmarks have been selected to include binary, categorical and ordinal decision variables but also different discrete structures such as sequential or graph-based structures. Since we assume that the evaluation of the objective functions is expensive, we perform the benchmark with a relatively strict budget of at most $500$ evaluations. The objective function is evaluated once per iteration in Algorithm~\ref{alg:surrogate_modelling}. 
Furthermore, we are testing on relatively large problem sizes, ranging from $44$ up to $150$ decision variables with search spaces of around $\sim10^{50}$ possibilities. This range is interesting considering that Bayesian optimization using Gaussian processes is typically applied on problems with less than 10 variables.

On top of that, it has been shown that a large dimensionality reduces the importance of choosing a complicated acquisition function~\cite{rehbach2020expected}, which helps us doing a fair comparison between surrogates.

Moreover, we do an analysis of the performance of each algorithm where we limit the allowed time budget instead of the number of evaluations and simulate different evaluation times of the objective functions. The time budget includes both the total time to evaluate the objective function and the computation time of the optimization algorithm. Thus, it puts emphasis on the computation time of the algorithm in addition to their respective sample efficiency.

We present the four benchmark problems in detail below. Note that we present these problems in detail but that they are treated as black boxes by the optimization algorithms.

\textbf{The Discrete Rosenbrock problem} is a $d$-dimensional, non-convex function, with a curved valley that contains the global optimum defined by the following function:
\begin{equation}
    f(\mathbf{x}) = \sum_{i=1}^{d-1} [100 (\mathbf{x}_{i+1} - \mathbf{x}_i^2 )^2 + (1-\mathbf{x}_i)^2]
\end{equation}
where $\mathbf{x}\in\mathbb{Z}^d$. In the Rosenbrock problem, finding the valley is simple, but finding the global optimum $[1,1,\ldots,1]$ is not. As we are exploring discrete optimization problems, we consider a discrete variant of the problem such that only integer solutions are considered. We have $d=49$ decision variables and each decision variable $x_i$ is bounded by the range $[-5,10]$. Thus, the problem's search space is in the order of $10^{59}$ candidate solutions. Lastly, the additive noise $\epsilon$ is normally distributed according to $N(\mu = 0,\sigma = 10
^{-6})$.  

\textbf{The Weighted Max-Cut problem} is an NP-hard graph cutting problem, defined as follows: For an undirected weighted graph $G = (V, E)$, a \emph{cut} in $G$ creates the subset $S \subseteq V$ and its complement $\overline{S} = V\backslash S$. Then $E(S,\overline{S})$, is defined as the set of edges that have one vertex in $S$ and the other in $\overline{S}$. The Max-Cut problem is to find the cut that maximizes the weight of the edges in $E(S,\overline{S})$. 
The problem is encoded with a binary string $x\in\{0,1\}^d$ where either $x_i=0$ or $x_i=1$ indicates if node $i$ lies in $S$ or $\overline{S}$ respectively.

For the following experiments, the MaxCut problem instances are randomly generated as weighted graphs, with $d$ nodes, edge probability $p=0.5$ and a uniformly distributed edge weight in the range $[0, 10]$.
The graph generator is initialized with the same random seed for every run, ensuring that all experiments of a given problem size are performed on the same graph.
On top of that, the additive noise $\epsilon$ added to each evaluation is following a standard normal distribution $N(\mu = 0, \sigma = 1)$. Lastly, we are using a graph with $d=150$ nodes which means that the size of the problem's search space is $2^{150}\approx 10^{57}$.

\textbf{The Perturbed Traveling Salesman} is a variant of the well-known sequential graph problem where, given a number of cities and the distances between these cities, a shortest path needs to be found that visits all cities and returns to the starting city. We consider the asymmetric case with $k$ cities where the distance between cities is not the same in both directions. Moreover, noise $\epsilon\sim U(0,1)$ is added to each distance during evaluation. While the perturbation can cause issues for problem-specific solvers, it creates a good benchmark for black-box optimization algorithms. To ensure a robust solution, each proposed route is also evaluated 100 times and the worst-case objective value is returned. 
Furthermore, we will consider problem instance \emph{ftv44}. This is an instance with 44 cities taken from TSPLIB \cite{TSPlib}, a library of problem instances for the traveling salesman problem. 
An instance with 44 cities is chosen to closely match the number of decision variables in the ESP problem which has a fixed number of 49 decision variables.

The problem is encoded as in~\cite{bliek2020black}: after choosing a fixed origin city, there are $d=k-2$ ordered decision variables $x_i$ for $i=1,\dots,d$ such that $x_1\in\{1,2,\dots,k-1\}$ where each integer represents a city other than the origin city. Then, the next decision variable $x_2\in\{1,...,k-2\}$ selects between the cities that were not yet visited. This is repeated until all cities have been chosen in some order. Since the last decision variable $x_{d}\in\{1,2\}$ selects between the two remaining cities, we can deduce afterward the two remaining edges which closes the route since there is one last city to visit before returning to the origin city. 
Thus, the total number of possible sequences is given by 
$(d-1)!\approx 6\cdot 10^{52}$ for this instance.

\textbf{The Electrostatic Precipitator problem} is a real-world industrial optimization problem first published by Rehbach et al.~\cite{rehbach2018ESPBenchmark}. The Electrostatic Precipitator (ESP) is a crucial component for gas cleaning systems. It is a large device that is used when solid particles need to be filtered from exhaust gases, such as reducing pollution in fossil fueled power plants. Before gas enters the ESP, it passes through a gas distribution system that controls the gas flow into the ESP. The gas flow is guided by configurable metal plates which blocks the airflow to a varying degree. The configuration of these plates inside the gas distribution system is vital for the efficiency of the ESP. However, it is non-trivial to configure this system optimally. 

The objective function is computed with a computationally intensive fluid dynamics simulation, taking about half a minute of computation time every time a configuration is tested. There are 49 slots where different types of plates can be placed or be left empty. In total, there are 8 different options available per slot. 
This is formalized such that each integer-valued solution $\mathbf{x}$ is subject to the inequality constraint $0\leq \mathbf{x}_i\leq 7$ for $i=1,\dots,49$. This gives a large solution space in the order of $10^{44}$ possibilities.  

Lastly, the problem has some ordinal structure where the decision variables decides between sizes of holes which are covering the plates. However, as an indication of the complex problem structure we have noted that changing any single variable does not affect the objective function.

\section{Experiments}
The goal of this section is to show a benchmark comparison between discrete and continuous surrogate-based algorithms on the discrete optimization problems of the previous section.
The compared algorithms are HyperOpt and SMAC as two popular surrogate-based algorithms that make use of a discrete surrogate model if the search space is discrete, and Bayesian optimization as a popular surrogate-based algorithm for continuous problems.
Though there exist several other algorithms that can deal with the discrete setting, these three are often used in practice because they are well established, can be used for a wide variety of problems, and have code available online.
The most recent method we found online, namely the PSR variant of BOCS, requires too much memory and computation time for problems of the size we consider in this work and is therefore not included in the comparison.
We do also include IDONE in the comparisons as a surrogate-based algorithm that uses a continuous surrogate model but is designed for discrete problems,
and random search is included as a baseline.

All experiments were run on the same Unix-based laptop with a Dual-Core Intel Core i5 2.7 GHz CPU and 8 GB RAM. Each algorithm attempted to solve the benchmarks $5$ times. The allowed evaluation budget was $500$ evaluations for all problems except the ESP problem where $100$ evaluations were allowed instead due to it being more computationally expensive. 

We are using the default hyperparameters for all algorithms, which are decided by their respective code libraries, with two exceptions. We change the SMAC algorithm to deterministic mode, since it otherwise evaluates the same point several times, which deteriorates its performance significantly.
Besides that, the first five iterations of IDONE are random evaluations, which is similar to what happens in the other algorithms.
The other algorithms start with their default number of random evaluations (which is $5$ for Bayesian optimization and $3$ for SMAC and HyperOpt), however for a fair comparison we make sure that all of these initial random evaluations come from a uniform distribution over the search space.
Unfortunately, more extensive hyperparameter tuning than stated above is  too time-consuming for expensive optimization problems such as ESP.

In the following section we present the results from the benchmark comparison of the four surrogate-based optimization algorithms. The benchmark consists of the four problems which have varying discrete structures.

\subsection{Results}

In this section we describe the main results from comparing the algorithms on the discrete Rosenbrock, weighted Max-Cut, the travelling salesman and the ESP problems. Figure~\ref{fig:benchmarks} shows the best average objective value found until a given iteration on each problem as well as their respective computation time. The computation time is the cumulative time up until iteration $i$ which is required to perform the steps on line 5 and 6 in Algorithm~\ref{alg:surrogate_modelling}. Furthermore, we also investigate how the algorithms perform if we introduce a time budget during optimization instead of constraining the number of evaluations. 

\subsubsection{Ordinal structures} 
\begin{figure}[tbp]
    \centering
    \begin{subfigure}{6cm}
        \centering\includegraphics[width=5.4cm]{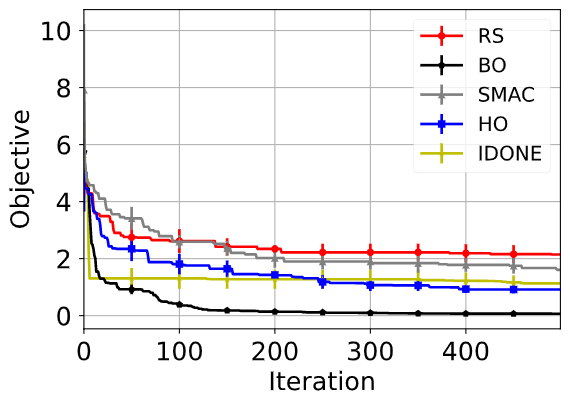}
        \caption{Best average objective value versus iteration on the 49-dimensional discrete Rosenbrock.}
        \label{fig:rosen49_benchmark}
    \end{subfigure}%
    \hspace{0.1cm}
    \begin{subfigure}{6cm}
        \centering\includegraphics[width=5.4cm]{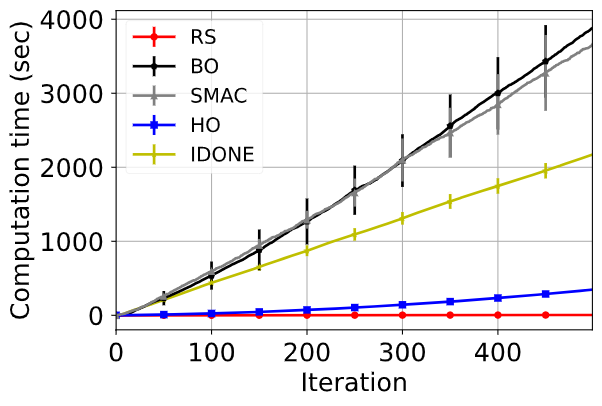}
        \caption{Average computation time versus iteration on the 49-dimensional discrete Rosenbrock.}
        \label{fig:rosen49_time}
    \end{subfigure}
    
    \begin{subfigure}{6cm}
        \centering\includegraphics[width=5.4cm]{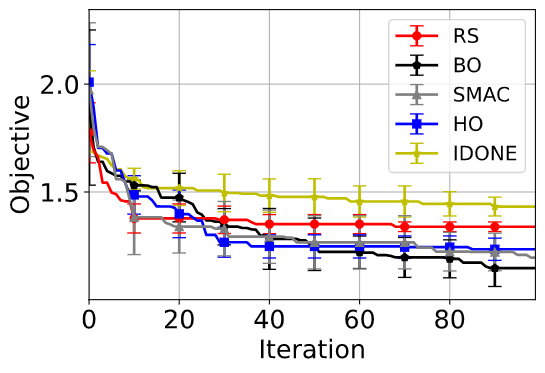}
        \caption{Best average objective value versus iteration on the ESP problem.}
        \label{fig:esp_benchmark}
    \end{subfigure}%
    \hspace{0.1cm}
    \begin{subfigure}{6cm}
        \centering\includegraphics[width=5.4cm]{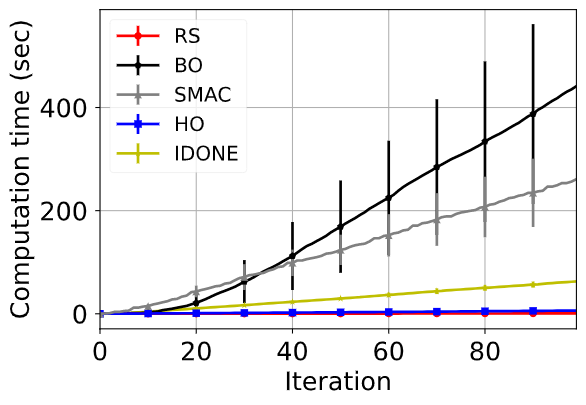}
        \caption{Average computation time versus iteration on the ESP problem.}
        \label{fig:esp_time}
    \end{subfigure}
    
    \begin{subfigure}{6cm}
        \centering\includegraphics[width=5.4cm]{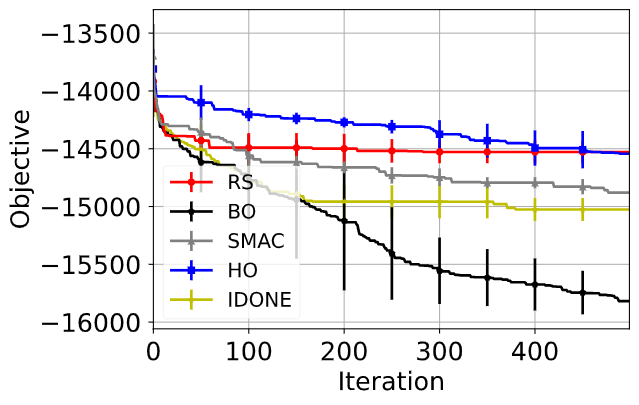}
        \caption{Best average objective value versus iteration on the 150-dimensional weighted Max-Cut.}
        \label{fig:maxcut150_benchmark}
    \end{subfigure}%
    \hspace{0.1cm}
    \begin{subfigure}{6cm}
        \centering\includegraphics[width=5.4cm]{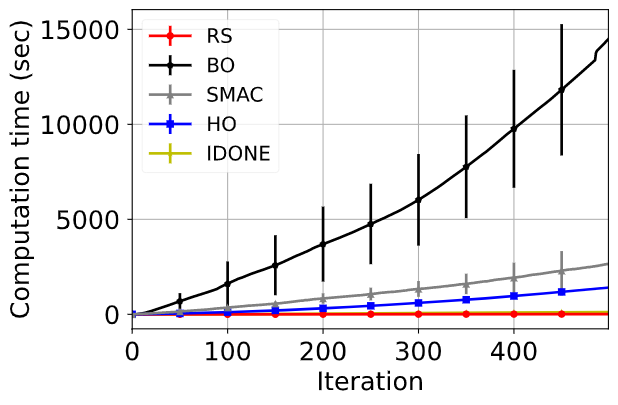}
        \caption{Average computation time versus iteration on the 150-dimensional weighted Max-Cut.}
        \label{fig:maxcut150_time}
    \end{subfigure}
    
    \begin{subfigure}{6cm}
        \centering\includegraphics[width=5.4cm]{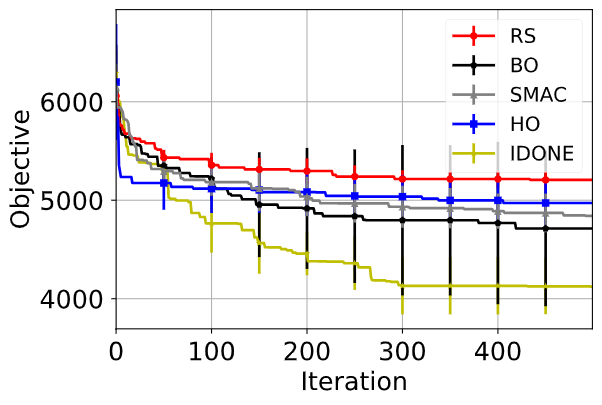}
        \caption{Best average objective value versus iteration on the TSP with 44 cities.}
        \label{fig:tsp44_benchmark}
    \end{subfigure}%
    \hspace{0.1cm}
    \begin{subfigure}{6cm}
        \centering\includegraphics[width=5.4cm]{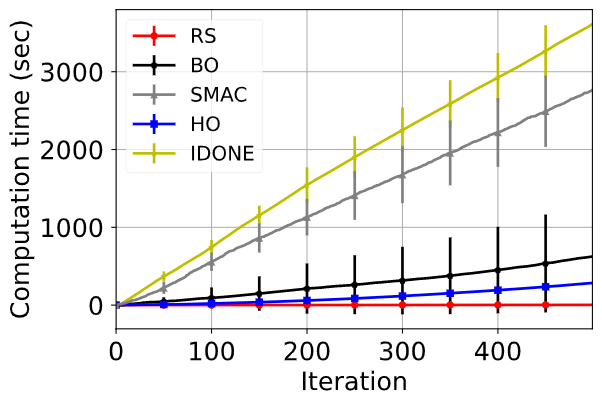}
        \caption{Average computation time versus iteration on the TSP with 44 cities.}
        \label{fig:tsp44_time}
    \end{subfigure}
    \caption{Comparison of objective value and computation time of Bayesian optimization (BO), SMAC, IDONE, HyperOpt (HO) and random search (RS) on four different benchmark problem.  An average is computed from 5 runs and the standard deviation is plotted as the error. The objective value has been negated for Max-Cut since the maximization problem has been turned into a minimization problem. The evaluation budget was 500 evaluations for all problems except the ESP problem which was limited to 100 evaluations due to it being more computationally expensive.}
    \label{fig:benchmarks}
\end{figure}

We start by comparing the results from the 49-dimensional discrete Rosenbrock problem.
In Figure~\ref{fig:rosen49_benchmark}, we see that Bayesian optimization (BO) is the only algorithm that comes close to the optimal objective value of zero. The other algorithms are not performing as well, where HyperOpt (HO) gets the closest to BO. Given that the problem is in fact a discrete version of an inherently continuous problem with ordinal variables, this can be considered to be well suited for continuous model regression. On the other hand, IDONE also uses a continuous surrogate, but it does not perform as well as BO. A possible explanation is that IDONE is less flexible since it is a piece-wise linear model. 

\begin{figure}
    \centering
    \includegraphics[width=0.9\textwidth]{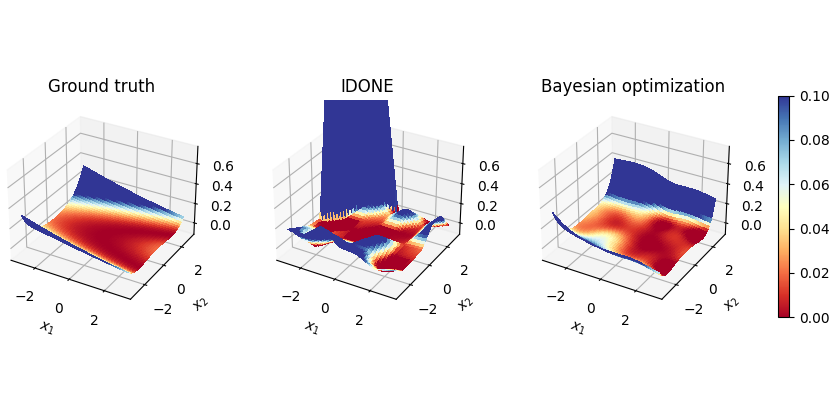}
    \caption{Visualization of continuous surrogates that approximate the two-dimensional Rosenbrock, namely the linear combination of ReLUs from IDONE and Gaussian processes from BO. These models were picked based on the best performance from 15 different runs with 50 evaluations each. HyperOpt and SMAC are not visualized since this is not supported by their respective code libraries.}
    \label{fig:surrogate_visualization}
\end{figure}

To investigate the quality of the surrogates from both BO and IDONE, we visualize their surfaces in Figure~\ref{fig:surrogate_visualization} for the $2$-dimensional case of Rosenbrock. The Gaussian process from BO (which uses a Mat\'ern 5/2 kernel in this case) predicts a smoother surface than IDONE which appears more rugged and uneven. Overall, BO looks more similar to the objective ground truth. We can argue that this is why BO performs well while IDONE does not. BO is likely suitable for the discrete Rosenbrock problem since the problem has an underlying continuous structure with ordinal variables. Meanwhile, this structure could be too complex for the piece-wise linear surrogate in IDONE.

However, we are interested in investigating problems which do not necessarily have a clear continuous structure. Thus, we look at the ESP problem which also happens to have some ordinal structure. The results from this problem are found in Figure~\ref{fig:esp_benchmark}. It shows a more even performance among the algorithms compared to the Rosenbrock problem, although BO still returns the best objective on average. This is closely followed by both SMAC and HO, while IDONE is doing worse than random search.

Based on the results from these two problems, it appears that BO works well on ordinal structures. However, this does not seem to hold true for all continuous surrogates considering the performance of IDONE. Still, the naive approach with BO outperforms the other state-of-the-art discrete algorithms on the problems that we have discussed so far.
This is actually in line with experimental results from~\cite{garrido2020dealing} on small problems (up to $6$ dimensions) with both discrete and continuous parameters, though it was not the main conclusion of the authors. 
The difference with our work is that we consider purely discrete problems of higher dimensions, from a real-life application, and we include IDONE in the comparison.

\begin{table}[tb]
\addtolength{\tabcolsep}{2pt}
\begin{center}
\begin{minipage}[c]{\textwidth}
    \begin{minipage}[t]{0.48\textwidth}
        \begin{table}[H]
            \centering
            \begin{tabular}{c|c|c}
                Algorithm & Non-binary & Binary \\ \hline
                BO      &  \textbf{0.067 (0.021)} &   \textbf{0.37 (0.038)}    \\
                SMAC    & 1.61 (0.18)  &  1.28 (0.29)     \\
                HyperOpt & 0.91 (0.13) &   0.94 (0.14)     \\
                IDONE   &  1.13 (0.20) &   0.61 (0.038)    \\
            \end{tabular}
            \vspace{4pt}
            \caption{Comparison of results on the 49-dimensional discrete Rosenbrock with and without binary encoding of the decision variables. The final average objective value from 5 runs is presented after 500 evaluations with the standard deviation in parenthesis. The lowest objective value is marked as bold in each column.}
            \label{tab:binarized_experiment_rosenbrock}
        \end{table}
    \end{minipage}
    \hfill
    \begin{minipage}[t]{0.48\textwidth}
        \begin{table}[H]
            \centering
   
            \begin{tabular}{c|c|c}
                Algorithm & Non-shuffled & Shuffled \\ \hline
                BO      &  4713.2 (789.2)  &  4898.0 (292.4)    \\
                SMAC    &  4841.8 (184.7) &  4784.9 (302.7)     \\
                HyperOpt &  4971.9 (256.5) &   4871.8 (221.9)   \\
                IDONE   & \textbf{4122.8 (279.8)} &   \textbf{4556.4 (175.7)}     \\
            \end{tabular}
            \vspace{4pt}
            \caption{Comparison of TSP with 44 cities when the input has a sequential structure versus that decision variables' position have been shuffled. The final average objective value from 5 runs is presented after 500 evaluations with the standard deviation in parenthesis. The lowest objective value is marked as bold in each column.}
            \label{tab:shuffled_TSP_experiment}
  
        \end{table}
    \end{minipage}
\end{minipage}
\end{center}
\addtolength{\tabcolsep}{-2pt}
\end{table}

\subsubsection{Binary structures}
We will now consider a graph problem, that is the weighted Max-Cut problem. From the results in Figure~\ref{fig:maxcut150_benchmark}, we notice that BO clearly outperforms all other algorithms. Meanwhile, IDONE is the second best, followed by SMAC and then HO which performs worse than random search. Compared to the other problems that we have seen so far, a major difference is the binary decision variables in the Max-Cut problem. We use this to frame our hypothesis, namely that the good performance of BO on the Max-Cut problem is due to the binary structure of the problem. 

To investigate this hypothesis, we perform an additional experiment by encoding the 49-dimensional, discrete Rosenbrock with binary variables and compare this with the previous results from Figure~\ref{fig:rosen49_benchmark}. The ordinary problem has 49 integer decision variables which lie in the range $[-5,10]$, this is converted into a total of 196 binary decision variables for the binary-encoded version. Table~\ref{tab:binarized_experiment_rosenbrock} shows the performance of the algorithms on the binary-encoded, discrete Rosenbrock. Although BO is performing worse on the binarized Rosenbrock, it is still performing the best compared to the other algorithms, even though both SMAC and IDONE perform better on the binarized problem.

Thus, we could argue that the binary representation of the Max-Cut problem can not explain why BO performs well on this problem. There is a possible argument that the binary variables might cause less rounding-off errors since the range of values is simply zero to one with a threshold in the middle. However, a counter-argument is that such a large number of decision variables is typically not well-suited for Gaussian processes regression. This is also indicated by the large computation time of BO on the Max-Cut, see Figure~\ref{fig:maxcut150_time}. 

\subsubsection{Sequential structures}
Even though TSP is a graph problem like the Max-Cut problem, there is an important difference. TSP has a sequential structure since the decision variables select an ordering that directly affects the objective value. Moreover, the encoding of the problem, as described in Section 4, causes strong interactions between adjacent decision variables. 

We continue by looking at the results from TSP in Figure~\ref{fig:tsp44_benchmark}. BO is now outperformed by IDONE even though it still performs better than SMAC and HO on average, although BO has a large variance on this problem. We suspect that the sequential structure is well-suited for IDONE, as it explicitly fits some of its basis functions with adjacent variables in the input vector $(x_1, x_2,\dots, x_d)$~\cite{bliek2020black}.

To investigate whether this is the case, we test what happens when the order of the decision variables are re-shuffled in TSP such that the sequential structure is removed.
This is done by adding to the objective function a mapping that changes the order of the variables in the input vector $(x_1, x_2,\dots, x_d)$ to a fixed arbitrarily chosen order.
From Table~\ref{tab:shuffled_TSP_experiment} we see that IDONE performs worse without the original sequential structure. At the same time, the other algorithms show no large significant difference.
However, IDONE returns the best objective on average both with and without shuffling the order of variables.
The large variance on BO makes it more difficult to draw any strong conclusions, but since IDONE also uses a continuous surrogate model, we can still conclude that continuous surrogates perform better than the discrete counterparts on this problem.

\subsubsection{Taking computation time into consideration}
Although BO performs well on the benchmark comparisons, we notice that it is more expensive with respect to computation time compared to the other surrogate-based methods. Figures~\ref{fig:rosen49_time}, \ref{fig:esp_time}, \ref{fig:maxcut150_time} and \ref{fig:tsp44_time} show the cumulative time on the problems.

In general, BO requires a vast amount of time compared to the other algorithms, especially on Max-Cut where the computations took one to two minutes per iteration. This is not surprising considering that regression with Gaussian processes is computationally intensive: its complexity grows as $O(n^3)$ where $n$ are the number of observations \cite{shahriari2016reviewBO}. This can be a big drawback if the evaluation time of the objective function is relatively small.

Meanwhile, the other algorithms share similar computation times which are often less than one second. The only exception is for IDONE which requires more computation time on TSP, see Figure~\ref{fig:tsp44_time}.

\begin{figure}[tb]
    \centering
    \includegraphics[width=\textwidth]{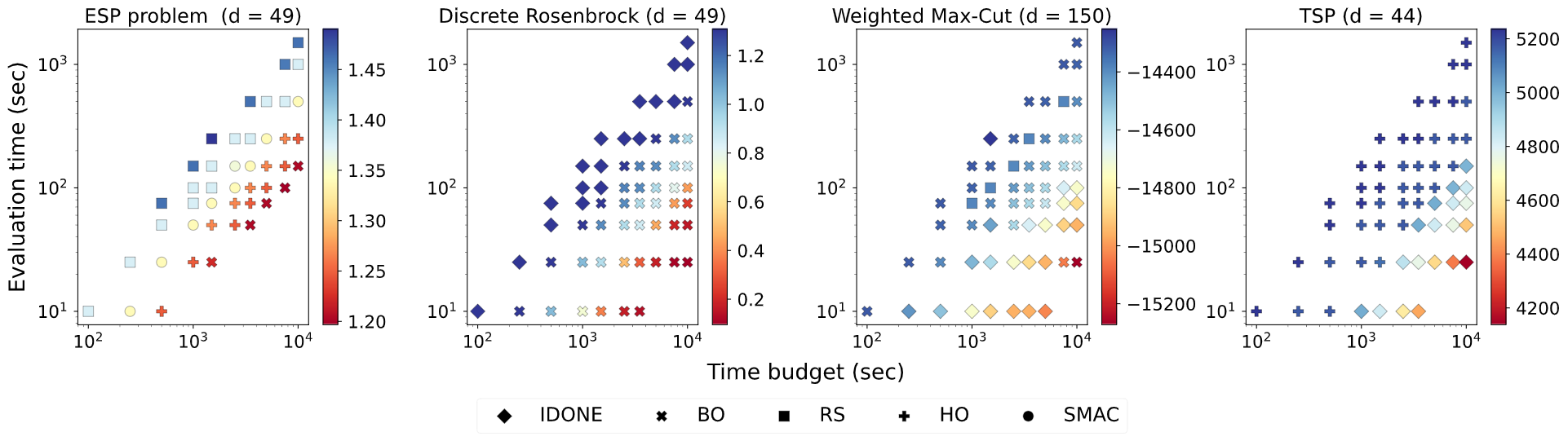}
    \caption{The best algorithm on average for a given time budget and an evaluation time of the objective function is indicated with different shapes for each algorithm. The colors represent the objective value. The time budget includes both the evaluation time and the computation time of the algorithms.
    These results are obtained by adding an artificial evaluation time after running the experiments.}
    \label{fig:eval_vs_time_budget}
\end{figure}

So far, we have only considered experiments that restrict the number of evaluations. But in real-life applications, the computation time of an algorithm can be important to take into consideration when limited with some given time budget as well. In particular, the large computation time of BO motivates the question whether it would still perform well under a constrained time budget instead.
By keeping track of both the evaluation times of the objective functions, as well as the computation time spent by the algorithms at every iteration, we can investigate the performance of the algorithms in different situations.
We artificially adjust the evaluation time in the experiments from Figure~\ref{fig:benchmarks} to simulate the cost of the objective function.
The evaluation time ranges from $10^1$ to $1.5\cdot 10^3$ seconds. Similarly, the time budget varies between $10^2$ and $10^4$ seconds.

Figure~\ref{fig:eval_vs_time_budget} displays which algorithm performs best on average for each problem, depending on the evaluation time and time budget. 
It also shows the objective value that was achieved by the best performing algorithm.
Results only occur below the line $y=x$ because the time budget must be larger than the evaluation time.
To ensure a fair comparison, we only present the algorithm with the best final average objective value if the maximum number of evaluations from the previous benchmark experiments was not exceeded within the allocated time budget for all algorithms.
As expected, on all problems, for a given evaluation time (value on the vertical axis), the objective values become lower (better) with an increase in the time budget.

For the ESP problem, the results are mixed. The best algorithm varies between BO, HO, SMAC and even random search, depending mostly on the ratio between the time budget and the evaluation time. For example, random search performs best when the evaluation time is around the order of $10^1$ smaller than the time budget which gives relatively few evaluations. Meanwhile, BO performs best with a much larger ratio.  On the discrete Rosenbrock benchmark, BO is clearly the best in almost all cases. The only exception is when the ratio between evaluation time and time budget is very small (so only 10-50 iterations can be performed), in which case IDONE performs better. For the weighted Max-Cut, on the other hand, we notice the opposite of what we see with the Rosenbrock benchmark. Thus, it seems like the growth in computation time of BO, see Figure~\ref{fig:maxcut150_time}, sometimes outweighs the good performance that we noted earlier when only taking an evaluation budget into consideration. Lastly, we see that IDONE and HyperOpt outperform other algorithms on TSP when constrained by a time budget. 

This experiment gives a better picture of the performance of each algorithm, especially if we may consider it to be more realistic by taking time constraints into consideration. Thus, the experiment from Figure~\ref{fig:eval_vs_time_budget} is a good complement to our benchmark comparison. 
In the following and last section, we summarize the conclusions that can be drawn from all of the above experiments. 

\section{Conclusion and Future Work}

Based on the results from the benchmark comparison, we can show that 
the use of continuous surrogate models is
a valid approach for expensive, discrete black-box optimization. Moreover, we give insight into what discrete problem structures are well-suited for 
the different methods.

We have shown that Bayesian optimization (BO) performs better than discrete state-of-the-art algorithms on the four tested, high-dimensional benchmarks problem with either ordinal, sequential or binary structures. IDONE, another continuous surrogate-based algorithm designed for discrete problems, outperforms BO on the benchmark with a sequential structure, but not on the three other benchmarks.

In addition, we have investigated how the different algorithms deal with the different problem structures. Firstly, ordinal structures appear suitable for BO, especially if the objective function has an underlying continuous structure such as the discrete Rosenbrock benchmark.
For binary structures, we noticed that BO is negatively affected by binary variables, while IDONE and SMAC benefited from this transformation. However, BO still returned the best solution on the binary Max-Cut problem, even though a big drawback was its computation time. Lastly, we have seen that 
IDONE outperforms the other algorithms on a problem with sequential decision variables, even after negatively affecting it by changing the ordering.

We also 
investigated the different algorithms under different time constraints by artificially changing the function evaluation times of the different benchmark problems.
For lower time budgets, BO is held back by its large computation time in some cases. 
Even though BO is a time-intensive method,
it mostly showed competitive performance when the evaluation time was relatively low and the time budget high, except for the binary Max-Cut problem.
IDONE, HyperOpt, SMAC, and even random search all had specific problems and time budgets where they outperformed other algorithms.
Lastly, based on our results, discrete surrogate-based methods could be more relevant in the setting with a limited time budget, in contrast to only limiting the number of evaluations.

Finally, we state some open questions which remain to be answered about continuous surrogates in the topic of expensive, discrete black-box optimization. 
Considering that we looked at a naive approach of BO, 
it is still an open question how the more advanced discrete BO variations would fare in the framework where time budgets and function evaluations times are taken into account like in this paper.
This same framework would also lead to interesting comparisons between surrogate-based algorithms and other black-box algorithms such as local search or evolutionary algorithms, which are better suited for cheap function evaluations.
It also remains unclear why BO performs
best on the binary Max-Cut benchmark even though it is negatively affected by binary structures on the Rosenbrock function.
Finally, it would be of great practical value if one could decide on the best surrogate-based algorithm in advance, given the time budget and evaluation time of a real-life optimization problem.
This research is a first step in that direction.

\section*{Acknowledgments}

This work is part of the research programme Real-time data-driven maintenance logistics with project number 628.009.012, which is financed by the Dutch Research Council (NWO).
The authors would also like to thank Arthur Guijt for helping with the python code.

\bibliographystyle{splncs04}
\bibliography{main}

\begin{thebibliography}{10}
\providecommand{\url}[1]{\texttt{#1}}
\providecommand{\urlprefix}{URL }
\providecommand{\doi}[1]{https://doi.org/#1}

\bibitem{baptista2018bayesian}
Baptista, R., Poloczek, M.: Bayesian optimization of combinatorial structures.
  In: International Conference on Machine Learning. pp. 471--480 (2018)

\bibitem{bergstra2013hyperopt}
Bergstra, J., Yamins, D., Cox, D.D.: Hyperopt: A {P}ython library for
  optimizing the hyperparameters of machine learning algorithms. In:
  Proceedings of the 12th Python in Science Conference. pp. 13--20 (2013)

\bibitem{bergstra2013making}
Bergstra, J., Yamins, D., Cox, D.D.: Making a science of model search:
  Hyperparameter optimization in hundreds of dimensions for vision
  architectures. In: Proceedings of the 30th International Conference on
  Machine Learning. Jmlr (2013)

\bibitem{bliek2020black}
Bliek, L., Verwer, S., de~Weerdt, M.: Black-box combinatorial optimization
  using models with integer-valued minima. Annals of Mathematics and Artificial
  Intelligence pp. 1--15 (2020)

\bibitem{de2019greed}
De~Ath, G., Everson, R.M., Rahat, A.A., Fieldsend, J.E.: Greed is good:
  Exploration and exploitation trade-offs in {B}ayesian optimisation. arXiv
  preprint arXiv:1911.12809  (2019)

\bibitem{deshwal2020scalable}
Deshwal, A., Belakaria, S., Doppa, J.R.: Scalable combinatorial {B}ayesian
  optimization with tractable statistical models. arXiv preprint
  arXiv:2008.08177  (2020)

\bibitem{elsken2018neural}
Elsken, T., Metzen, J.H., Hutter, F.: Neural architecture search: A survey.
  arXiv preprint arXiv:1808.05377  (2018)

\bibitem{garrido2020dealing}
Garrido-Merch{\'a}n, E.C., Hern{\'a}ndez-Lobato, D.: Dealing with categorical
  and integer-valued variables in {B}ayesian optimization with {G}aussian
  processes. Neurocomputing  \textbf{380},  20--35 (2020)

\bibitem{griffiths2017constrained}
Griffiths, R.R., Hern{\'a}ndez-Lobato, J.M.: Constrained {B}ayesian
  optimization for automatic chemical design. arXiv preprint arXiv:1709.05501
  (2017)

\bibitem{hernandez-lobato2017NIPSworkshop}
Hern{\'a}ndez-Lobato, J.M., Gonzalez, J., Martinez-Cantin, R.: {NIPS} workshop
  on {B}ayesian optimization. \url{https://bayesopt.github.io/}, accessed
  22-08-2020

\bibitem{hernandez2014predictive}
Hern{\'a}ndez-Lobato, J.M., Hoffman, M.W., Ghahramani, Z.: Predictive entropy
  search for efficient global optimization of black-box functions. In: Advances
  in Neural Information Processing Systems. pp. 918--926 (2014)

\bibitem{hutter2011sequential}
Hutter, F., Hoos, H.H., Leyton-Brown, K.: Sequential model-based optimization
  for general algorithm configuration. In: International Conference on Learning
  and Intelligent Optimization. pp. 507--523. Springer (2011)

\bibitem{jones1998efficient}
Jones, D.R., Schonlau, M., Welch, W.J.: Efficient global optimization of
  expensive black-box functions. Journal of Global optimization
  \textbf{13}(4),  455--492 (1998)

\bibitem{klein2017fast}
Klein, A., Falkner, S., Bartels, S., Hennig, P., Hutter, F.: Fast {B}ayesian
  optimization of machine learning hyperparameters on large datasets. In:
  Artificial Intelligence and Statistics. pp. 528--536 (2017)

\bibitem{korovina2020chembo}
Korovina, K., Xu, S., Kandasamy, K., Neiswanger, W., Poczos, B., Schneider, J.,
  Xing, E.: Chembo: Bayesian optimization of small organic molecules with
  synthesizable recommendations. In: International Conference on Artificial
  Intelligence and Statistics. pp. 3393--3403. PMLR (2020)

\bibitem{lindauer2017warmstarting}
Lindauer, M., Hutter, F.: Warmstarting of model-based algorithm configuration.
  arXiv preprint arXiv:1709.04636  (2017)

\bibitem{luong2019bayesian}
Luong, P., Gupta, S., Nguyen, D., Rana, S., Venkatesh, S.: {B}ayesian
  optimization with discrete variables. In: Australasian Joint Conference on
  Artificial Intelligence. pp. 473--484. Springer (2019)

\bibitem{oh2019combinatorial}
Oh, C., Tomczak, J., Gavves, E., Welling, M.: Combinatorial bayesian
  optimization using the graph cartesian product. In: Advances in Neural
  Information Processing Systems. pp. 2914--2924 (2019)

\bibitem{rehbach2020expected}
Rehbach, F., Zaefferer, M., Naujoks, B., Bartz-Beielstein, T.: Expected
  improvement versus predicted value in surrogate-based optimization. arXiv
  preprint arXiv:2001.02957  (2020)

\bibitem{rehbach2018ESPBenchmark}
Rehbach, F., Zaefferer, M., Stork, J., Bartz-Beielstein, T.: Comparison of
  parallel surrogate-assisted optimization approaches. In: Proceedings of the
  Genetic and Evolutionary Computation Conference. p. 1348–1355. GECCO ’18,
  Association for Computing Machinery (2018)

\bibitem{TSPlib}
Reinelt, G.: {TSPlib}.
  \url{http://elib.zib.de/pub/mp-testdata/tsp/tsplib/tsplib.html}, accessed
  31-07-2020

\bibitem{shahriari2016reviewBO}
{Shahriari}, B., {Swersky}, K., {Wang}, Z., {Adams}, R.P., {de Freitas}, N.:
  Taking the human out of the loop: A review of {B}ayesian optimization.
  Proceedings of the IEEE  \textbf{104}(1),  148--175 (2016)

\bibitem{snoek2012practical}
Snoek, J., Larochelle, H., Adams, R.P.: Practical {B}ayesian optimization of
  machine learning algorithms. In: Advances in Neural Information Processing
  Systems. pp. 2951--2959 (2012)

\bibitem{ueno2016combo}
Ueno, T., Rhone, T.D., Hou, Z., Mizoguchi, T., Tsuda, K.: Combo: An efficient
  {B}ayesian optimization library for materials science. Materials discovery
  \textbf{4},  18--21 (2016)

\bibitem{zaefferer2018Auth}
Zaefferer, M.: Surrogate Models For Discrete Optimization Problems. Ph.D.
  thesis, Technische Universität Dortmund (2018)

\end{thebibliography}

\end{document}